\documentstyle[12pt]{article}
 \parskip
0.2cm 
\begin{document}

\newtheorem{guess}{Theorem}[section]
\newcommand{\bth}{\begin{guess}$\!\!\!${\bf .}~}
\newcommand{\eth}{\end{guess}}

\newtheorem{propo}[guess]{Proposition}
\newcommand{\bprop}{\begin{propo}$\!\!\!${\bf .}~}
\newcommand{\eprop}{\end{propo}}

\newtheorem{lema}[guess]{Lemma}
\newcommand{\blem}{\begin{lema}$\!\!\!${\bf .}~}
\newcommand{\elem}{\end{lema}}

\newtheorem{defe}[guess]{Definition}
\newcommand{\bdefe}{\begin{defe}$\!\!\!${\bf .}~}
\newcommand{\edefe}{\end{defe}}

\newtheorem{coro}[guess]{Corollary}
\newcommand{\bcor}{\begin{coro}$\!\!\!${\bf .}~}
\newcommand{\ecor}{\end{coro}}

\newtheorem{rema}[guess]{Remark}
\newcommand{\brem}{\begin{rema}$\!\!\!${\bf .}~\rm}
\newcommand{\erem}{\end{rema}}

\newcommand{\epf}{\hfill$\diamondsuit$}

\newcommand{\sing}{{\rm Sing}\,}
\newcommand{\lr}{\longrightarrow}
\newcommand{\ce}{{\cal E}}
\newcommand{\co}{{\cal O}} \newcommand{\qa}{{\cal Q}}
\newcommand{\ct}{{\cal T}}
\newcommand{\cf}{{\cal F}} \newcommand{\cg}{{\cal G}}
\newcommand{\cv}{{\cal V}}
\newcommand{\li}{{\cal L}} \newcommand{\cn}{{\cal N}}
\newcommand{\cu}{{\cal U}}
\newcommand{\ch}{{\cal H}} \newcommand{\cm}{{\cal M}}
\newcommand{\cc}{{\cal C}}
\newcommand{\cw}{{\cal W}} \newcommand{\cp}{{\cal P}}
\newcommand{\cl}{{\cal L}} \newcommand{\fq}{{q}}
\newcommand{\cd}{{\cal D}}
\newcommand{\Sing}{{\rm Sing}\,}
\newcommand{\eee}{{\rm End}\,}
\newcommand{\rank}{{\rm rank}\,}
\newcommand{\codim}{{\rm codim}\,}
\newcommand{\coker}{{\rm coker}\,}
\newcommand{\Aut}{{\rm Aut}\,}
\newcommand{\Ext}{{\rm Ext}\,}
\newcommand{\ad}{{\rm ad}\,}
\newcommand{\lm}{{\mbox{\Large $|$}}}
\newcommand{\rk}{{\rm rk}\,}
\newcommand{\wn}{{\cal W}^{k-1}_{n,d}}
\newcommand{\mo}{{\cal M}(n,d)}
\newcommand{\wl}{{\cal W}^{k-1}_{1,d}}
\newcommand{\ok}{{\cal O}^k}
\newcommand{\ra}{\rightarrow}
\newcommand{\da}{\downarrow}
\newcommand{\wnn}{{\cal W}^{k}_{n,d}}
\newcommand{\wnt}{\widetilde{{\cal
W}}^{k-1}_{n,d}} \newcommand{\wt}{\widetilde}
\newcommand{\mot}{\widetilde{{\cal
M}}(n,d)} \newcommand{\rn}{\rho^{k-1}_{n,d}}
\newcommand{\rl}{\rho^{k-1}_{1,d}}
\newcommand{\ctext}[1]{\makebox(0,0){#1}}
\newcommand{\Bbb}{\bf}
\setlength{\unitlength}{0.1mm}

\title{ GEOGRAPHY OF BRILL-NOETHER LOCI FOR SMALL SLOPES}
\author{L.
Brambila-Paz, I. Grzegorczyk, P. E. Newstead\thanks{All
the authors are members of the  VBAC Research group of
Europroj.  The work forms part of a project ``Moduli of
bundles in algebraic geometry'', funded by the EU
International Scientific  Cooperation Initiative (Contract
no. CI1$^*$-CT93-0031).  The first author also
acknowledges support from CONACYT (Project no. 3231-E9307)
and the second author  from UMD Foundation and EPSRC
(Grant no. GR/K76665).}}\date{October 1995}
\maketitle

\begin{abstract} {\small Let $X$ be
a non-singular  projective curve of genus
$g\ge2$ over an algebraically closed field of
characteristic zero. Let $\mo$
denote the moduli space of stable bundles of rank $n$ and
degree $d$ on $X$  and
$\wn $ the Brill-Noether loci in $\mo .$ We prove that,
if $0\leq d \leq n $ and  $\wn $ is non-empty, then it is
irreducible of the expected dimension and smooth  outside
$\wnn$. We prove further that in this range $\wn$ is
non-empty if and  only  if $d>0$, $n\leq d+(n-k)g$ and
$(n,d,k) \not= (n,n,n)$. We also prove   irreducibility
and non-emptiness for the semistable Brill-Noether loci.}
\end{abstract} \footnotetext[0]{1991 Mathematics Subject
Classification.  Primary: 14J60. Secondary: 14D20. }
\pagebreak \section*{Introduction}

\bigskip The moduli spaces of stable vector bundles over
an algebraic curve  have been extensively studied from
many points of view since they were first  constructed
more than 30 years ago, and much is now known about their
detailed  structure, in particular in terms of their
topology. Except in the classical  case of line bundles,
however, relatively little is known about their geometry
in terms, for example, of the existence and structure of
their subvarieties.

In the case of line bundles, where the moduli spaces are
all isomorphic to the  Jacobian, Brill-Noether theory has
long provided a basic source of geometrical  information.
This theory, which originated in the last century, is
concerned  with the subvarieties of the moduli spaces
determined by bundles having at least  a specified number
of independent sections. Basic questions, concerning
non-emptiness, connectedness, irreducibility, dimension,
singularities, cohomology  classes, etc., have been
completely answered when the underlying curve is  generic,
and departures from the generic behaviour are indeed used
to describe  curves with special properties.

The definitions can easily be extended to bundles of any
rank, but the basic  questions are then far from being
answered even for a generic curve. In  particular, given
integers $n$, $d$, $k$ with $n\geq2$ and $k\geq1$, one
would  like to know when there exist stable (or
semistable) bundles of rank $n$ and  degree $d$ having at
least $k$ independent sections. We use the term \lq\lq
geography'' to refer to the study of this problem (and
related questions such as  irreducibility and dimension of
the corresponding loci) by analogy with a  similar use of
the term in the theory of algebraic surfaces; we shall see
indeed  that much of the data obtained by ourselves and
others can be conveniently  summarised in graphical form
(see \S2 and particularly Figures 1 and 2).

In order to describe some of these ideas in more detail
and to state our own  results, we introduce some further
notation. Let $X$ be a non-singular projective  curve of
genus $ g\geq 2 $ defined over an algebraically closed
field of  characteristic $0$ and let $\mo $ denote the
moduli space of stable vector  bundles over $X$ of rank
$n$ and degree $d$. For any integer $k\geq1$, the {\it
Brill-Noether locus} $\wn$ is the set of stable bundles in
$\mo$ having at least  $k$ independent global sections;
this is in fact a subvariety of $\mo$ (see  1.1). (The
superscript $k-1$ is used here rather than $k$ largely for
historical reasons, since projective dimension rather than
vector space dimension was regarded classically as the
important notion.) Associated with this locus is the
number $$\rn = n^2(g-1)+1-k(k-d+n(g- 1)).$$ This is called
the {\it Brill-Noether number}, and is the \lq\lq
expected'' dimension of $\wn$. In a similar way, we denote
by $\mot$ the moduli  space of S-equivalence classes of
semistable vector bundles over $X$ and by   ${\wnt} $ the
corresponding Brill-Noether locus; again this is a
subvariety of  $\mot$ (see 1.2 for full definitions in
this case).

In the case $n=1$, the Brill-Noether loci (as remarked
above) have been well  known since the last century. In
fact the variety $\wl$ is always non-empty if  $\rl\geq0$,
and connected if $\rl>0$. The variety may be  reducible,
but each  component has dimension at least $\rl $. For a
generic curve $X$, $\wl$ is empty  if $\rl<0$ and is
irreducible of dimension $\rl$ with singular locus ${\cal
W}^{k}_{1,d}$ if  $g>\rl>0$. Modern proofs of these results
have been given by Kempf, Kleiman and  Laksov, Fulton and
Lazarsfeld, Griffiths and Harris, and Gieseker. A full
treatment of this case is contained in [ACGH].

For higher rank and $X$ generic, it is known that, for
$0<d\leq n(g-1)$, ${\cal  W}^0_{n,d}$ is irreducible of
dimension $\rho^0_{n,d}$ [Su] and $\Sing{\cal
W}^0_{n,d}={\cal W}^1_{n,d}$ [L]. The most extensive
results to date are  those of Teixidor [Te2]; these
describe many  cases when $\rn \geq 0$ and $\wn$ is
non-empty, as expected (see 2.5, where we shall use these
results  to draw our \lq\lq map''). Further results on
non-emptiness and irreducibility  are known when $n=2$ and
$k=2,3$ [Su, T, Te1, Te3], while ${\cal W}^{k-1}_{3,1}$
and  ${\cal W}^{k-1}_{3,2}$ are described in [NB]. On the
other hand, even for $X$  generic, $\wn$ may have
components of dimension greater than $\rn$ [BF] and the
singular set of $\wn$ may be strictly larger
than $\wnn$ [Te2].

In this paper we consider the case when  $n\geq 2$ and $0\leq d\leq n$ and
study
the varieties $\wn $ and $\wnt$. Our main results, which provide a complete
answer to the basic questions in  this case, are:

\bigskip

{\bf THEOREM A :} {\it If $\wn$ is non-empty, then it is
irreducible, of   dimension $\rn $ and $\Sing \wn ={\cal
W}^k_{n,d} $.}

\bigskip

{\bf THEOREM $\tilde{{\rm {\bf
A}}}$ :} {\it If $\wnt$ is non-empty, then it is
irreducible.}

\bigskip

{\bf THEOREM B :} {\it $\wn$ is non-empty if and
only if $$d>0,\ \  n\leq  d+(n-k)g\ \ and\ \
(n,d,k)\neq(n,n,n).$$ }

{\bf THEOREM $\tilde{{\rm {\bf B}}}$:} {\it  $\wnt $ is
non-empty  if and only  if either $$d=0\ \ and\ \  k\leq
n$$or $$d>0\ \  and\ \   n\leq d+(n-k)g.$$}

Our results give partial answers to questions 1 and 3 on
the VBAC Problems List  [VBAC], and are valid for all
non-singular curves, not just generic ones. Note  that, in
the case $k\leq d<n$, Theorems B and $\wt{{\rm  B}}$
follow from Teixidor's  results [Te2]. Note also that the
condition $n\leq d+(n-k)g $ implies that $\rn  \geq 1,$ so
in particular $\wn $ is empty when $\rn =0$ for $ 0\leq d
\leq n$;  this gives another example where the results in
higher rank differ from those in  rank 1. After the work
for this paper was completed, an alternative proof of
Theorem $\wt{{\rm  B}}$, using variational methods based
on the Yang-Mills-Higgs functional, was announced by G.
Daskalopoulos and R.  Wentworth [DW].

In proving our theorems, we shall distinguish the three
cases $0<d<n$, $d=0$ and  $d=n$, although all three will
depend on the use of extensions of the form $$0\rightarrow
{{\cal O}^k}\rightarrow E \rightarrow F \rightarrow 0.$$
In \S 1 we fix notation and give the basic definitions. In
\S 2 we give a proof (due to G. Xiao) of Clifford's Theorem
for semistable bundles (Theorem 2.1) and explain the
geography of the Brill-Noether loci. In \S3, we introduce
the use of extensions  (Proposition 3.1) and prove the
necessary condition  $n\leq d+(n-k)g$ in  Theorems B and
$\wt{{\rm  B}}$ (Theorem 3.3). In  \S 4 we prove Theorems
A and   $\wt{{\rm A}}$ when $0<d<n$ (Theorems 4.3, 4.4).
\S5 provides the setting for  the  proofs of Theorems B and
$\wt{{\rm  B}}$ which are completed in \S 6  (Theorem
6.3).  Finally, in \S\S 7, 8,  we prove all four theorems
for $d = 0 $  (Theorems 7.1, 7.3) and $d=n$ (Theorems 8.2,
8.5).

Our methods yield some information on the more detailed
geometry of the Brill- Noether loci (see, for example,
Corollary 4.5, Theorem 7.2, Theorem 8.3). These  varieties
are also closely connected with various types of augmented
bundle for  which moduli spaces have recently been
constructed. These include $k$-pairs [BeDW], coherent
systems [LeP1, 2] (also discussed as \lq\lq Brill-Noether
pairs'' in  [KN], and just \lq\lq pairs'' in [Be, RV]) and
extensions [BG]; for a general  survey, see [BDGW]. We
propose to return to these questions in future papers.

{\bf Acknowledgement.} The work for this paper was
completed during a visit by  the first two authors to
Liverpool. They wish to acknowledge the generous
hospitality of the University of Liverpool. The second
author would like to thank Institut Henri Poincar\'{e}
for support and hospitality. All the authors wish to
thank A.~D.  King for many  useful discussions; the
introductory material in \S\S1, 2 (and in  particular the
map of \S2) owe a great deal to unpublished notes of King.
We  also wish to thank R. Morris for designing Figures 1
and 2, and M.  Tapia (CIMAT) for help with
computer calculations which refuted an earlier  conjecture
and helped to lead us to a correct statement of
Proposition 6.1.

\bigskip

\renewcommand{\thesection}{\S \arabic{section}}
\section{Notation and
definitions}\renewcommand{\thesection}{\arabic{section}}

In this section, we give some basic notations and
definitions.

We denote by $X$ a non-singular projective curve of genus
$g\geq 2$,  fixed  throughout the paper, and write ${\ok
}={\cal O }^k_X$ for the trivial bundle of  rank $k$ over
$X$. For any integers $n$ and $d$ with $n\geq1$, let
${\mo}$  denote the moduli space of stable vector  bundles
of rank $n$ and degree $d$  over $X$. We write
$\mu(E)=\deg E/\rk E$ for the {\it slope} of a bundle $E$.

We make no distinction between locally free sheaves and
vector bundles over $X$.  However a subsheaf of a vector
bundle is called a subbundle only if the quotient  is
itself a vector bundle.

\bigskip  {\bf 1.1.} {\it Brill-Noether loci for stable
bundles. } As a set of  points, ${\wn}$ can be defined by
$$ {\wn}=\{ E\in {\mo}| h^0(E)\geq k\} .$$

Suppose first $(n,d)=1$. To obtain a scheme structure on
${\wn}$, let  ${\cal  U}$ be a universal bundle over $
X\times {\mo}$. Choose an effective divisor $D$  of  a
sufficiently large degree that $H^1(E\otimes L(D))=0$ for
all $E\in  {\mo}$. (Here $L(D)$ is the line bundle
associated to $D$.) Then,  in the exact  sequence
$$0{\ra}H^0(E){\ra}H^0(E\otimes
L(D)){\ra}H^0(E|_D){\ra}H^1(E){\ra}0,$$  the middle two
terms have dimensions  independent of $E$.
Globalising this, we obtain
$$0\ra\pi_*\cu\ra\pi_*(\cu\otimes
p^*_XL(D))\stackrel{\phi}{\ra}
\pi_*\left(\cu|_{D\times\mo}\right)\ra R^1_{\pi}\cu\ra
0,$$where $\pi:X\times  {\mo}\ra\mo$ is the projection
map. The middle two terms of this sequence are vector
bundles.

We can now define $\wn$ as the determinantal locus where
$\phi$ drops rank by at  least $k$. The \lq\lq expected''
dimension of $\wn$ is given by
$$\rn=n^2(g-1)+1-k(k-d+n(g-1)),$$ which is called the
{\it Brill-Noether number} associated to $\wn$. It follows
from the theory of determinantal varieties (see [ACGH] for
further details)  that, if $\wn$ is non-empty and
$\wn\neq\mo$, then $\dim\wn\geq\rn$.

For $(n,d)\not=1$, there is no universal bundle over
$X\times {\mo}$. However  the above construction works for
any locally universal family (for instance,  over a Quot
scheme); we can then define $\wn$ to be the image of the
variety so  obtained under the natural morphism to $\mo$.
It follows from geometric  invariant theory that this is a
closed subvariety of $\mo$.  \bigskip

{\bf 1.2. }{\it Brill-Noether loci for semistable bundles.}

Let $E$ be a semistable bundle of rank $n$. Then there
exists a filtration
$$0=E_0\subset E_1\subset E_2......\subset E_r=E$$ such
that $E_i/E_{i- 1}$ is a stable bundle with $\mu
(E_i/ E_{i-1})= \mu (E)$  for $0<i\leq r$.   The
associated graded bundle  $\bigoplus_i (E_i/ E_{i-1})$
depends only on $E$  and is denoted by ${\rm gr}\, E$.

We say that two semistable  bundles are {\it S-equivalent} if ${\rm gr}\,
E\cong{\rm gr}\, F$. There exists a moduli space
$\widetilde {\cal M}(n,d)$ of  S-equivalence classes of
semistable vector  bundles of rank $n$ and degree $d$,
which is an irreducible projective variety and is a
natural compactification of  $\mo$ (see [S]).

Writing $[E]$ for the S-equivalence class of $E$, we now
define  $$ \wnt=\{  [E]\in {\mo}| h^0({\rm gr}\, E)\geq
k\} .$$ Since $h^0({\rm gr}\, E)\geq h^0(E)$  for all $E$,
we can also define $\wnt$ as the set of S-equivalence
classes which  contain a bundle $E$ with $h^0(E)\geq k$.
We can give $\wnt$ a structure of  variety by using a
locally universal family as above. Note that this variety
does not have to be the closure of  ${\wn}$, as there may
exist components  containing semi-stable bundles only.
These components may have dimension smaller  then $\rn$.
For examples where this occurs, see \S7.

\bigskip {\bf 1.3. }{\it Petri map}. If $h^0(E)=k$, the
tangent space to ${\wn}$  at $E$ is the kernel of the map
$$ p^*: {\rm Ext^1} (E,E){\ra }H^0(E)^*\otimes  H^1(E)$$
which is dual  to the {\it Petri map} $$ p: H^0(E)\otimes
H^0(E^*\otimes K){\ra } H^0({\rm End}(E)\otimes K)$$
defined by multiplication  of sections. It follows easily
that ${\wn}$ is smooth of dimension $\rn$ at $E$  if and
only if the Petri map is injective. (Incidentally there
exist bundles $E$  for which the Petri map is not
injective [Te2, \S5].)

Note also that $\wnn\subset\Sing\wn$ whenever
$\wn\neq\mo$ (see [ACGH, Chapter II \S2 and p.~189).

 \bigskip
\renewcommand{\thesection}{\S\arabic{section}}\section{
\bf Brill-Noether geography of vector bundles of higher
ranks} \renewcommand{\thesection}{\arabic{section}} Our
main object in this section is to produce a \lq\lq map''
on which we can  display the results of Brill-Noether
theory for bundles of arbitrary rank.  Before doing this,
however, we shall state and prove a simple but
fundamental  result, which is a direct generalisation of
Clifford's Theorem for line bundles.

\noindent{\bf Theorem 2.1 (Clifford's Theorem). } {\it
Let $E$ be a semistable  bundle of rank $n$ and degree $d$
with $0\leq\mu(E)\leq 2g-2$. Then $$h^0(E)\leq
n+{d\over2}.$$}

{\it Proof:} (As far as we are aware, no complete proof of
this result has  appeared in the literature. The following
is due to G. Xiao.)

The proof is by induction on $n$, the case $n=1$ being the classical theorem.
For $E$ a semistable bundle of rank $n\geq2$, note first
that we can assume that  $h^0(E)>0$ and $h^1(E)>0$ (i.e.
$E$ is {\it special}), for otherwise the result  follows
at once from Riemann-Roch.

Now let $E_1$ be a proper subbundle of $E$ of maximal
slope and let $E_2=E/E_1$.  Certainly $E_1$ and $E_2$ are
both semistable. By semistability of $E$, we have
$\mu(E_1)\leq 2g-2$ and $\mu(E_2)\geq0$. On the other
hand, since $h^0(E)>0$,  $E$ possesses a subbundle of
non-negative degree; so $\mu(E_1)\geq0$. Similarly,  since
$h^1(E)>0$, $E$ possesses a quotient line bundle of degree
$\leq 2g-2$; by  comparing the slope of the kernel of this
quotient with that of $E_1$, one sees  easily that
$\mu(E_2)\leq 2g-2$. The result now follows at once by
induction.\epf

 \bigskip To construct our map, we first associate with
$\wn$ and $\wnt$ the  rational numbers $$\lambda={k\over
n},\ \ \mu={d\over n}.$$If $d<0$ and $k>0$, then $\wnt$ is
empty, while obviously $\wn=\mo$ if $k\leq 0$. We can
therefore plot $\mu$  against $\lambda$ in  the first
quadrant of the standard coordinate system (see  Figure
1). Advantages of plotting things in this way are that
every point with  rational coordinates can in principle
support bundles and that all ranks are  represented in the
same diagram.

In the remainder of the section, we describe some
important features of the map.

\bigskip {\bf 2.2}  {\it Riemann-Roch line} $\mu=
\lambda+g-1$.  By the Riemann- Roch Theorem
$$h^0(E)-h^1(E)=d-n(g-1).$$ Therefore for $\mu \geq
\lambda +g-1$,  i.e. above the Riemann-Roch line,
${\wn}$  is the whole moduli space.    Note  also that any
semistable bundle $E$ with $\mu(E)>2g-2$ has $h^1(E)=0$;
so, for  $\mu>2g-2$, $\wnt$ is empty below the
Riemann-Roch line.

\bigskip {\bf 2.3.} {\it Clifford line} $\mu=2\lambda-2$.
By Theorem 2.1  every  ${\wn}$  below this line is empty.

The interesting part of the map is therefore the
pentagonal region bounded by  the axes, the Riemann-Roch
line, the Clifford line and the line $\mu=2g-2$. This
corresponds to the region in which there may exist special
semistable bundles.

\bigskip{\bf 2.4.} {\it Brill-Noether curve}. Define
$$\tilde{\rho} ={1\over n^2}(\rn-
1)=(g-1)-\lambda(\lambda-\mu+(g-1)).$$We call the curve
$\tilde{\rho}=0$ the  {\it Brill-Noether curve}. The curve
is a branch of a hyperbola, below which the  expectation
is that the Brill-Noether loci will be finite.

\bigskip {\bf 2.5.} {\it Teixidor parallelograms}. In
[Te2], Teixidor defines  ranges of values for $n$, $d$,
$k$ such that for generic curves the  ${\wn}$ are
non-empty and have a component of the expected dimension.
These ranges  correspond to points $(\lambda, \mu )$ lying
in or on one of the parallelograms  marked T on the map.
These parallelograms have vertices at integer points,
sides  parallel to $\lambda=0$ and $\mu=\lambda$ and have
all their vertices on or  above the Brill-Noether curve
2.4. If all the vertices lie above      $\tilde{\rho} =
0$, then $\wn$ is non-empty whenever $(\lambda, \mu )$
lies in  or on the parallelogram. If the lower right
vertex of the parallelogram lies on  $\tilde{\rho} = 0$,
this still holds with the possible exception of those
points  of the parallelogram with the same
$\mu$-coordinate as this vertex; for such  points,
Teixidor shows only that $\wnt$ is non-empty.

\bigskip\bigskip\centerline{FIGURE 1}

\bigskip\bigskip In this paper we are concerned with the
region   $ {0\leq \mu \leq 1}$  of the map (see Figure 2).
The subregion $\lambda\leq\mu<1$ lies in a Teixidor
parallelogram, but the remainder of the region does not.
In any case, Teixidor  proves only non-emptiness (and, for
$X$ generic, the existence of a component of  the correct
dimension), whereas we shall solve the non-emptiness,
irreducibility  and singularity problems for the entire
region.

A key r\^{o}le in this is played by  the tangent line at
$(1,1)$ to   $\tilde{\rho}=0$. This is given by
$\mu+(1-\lambda)g=1$ or equivalently  $n=d+(n-k)g$. Thus
the inequality $n\leq d+(n-k)g$ in Theorems B and
$\wt{{\rm   B}}$ describes the area on or above this
tangent line. Theorem B therefore  states that for $ \mu
\leq 1$, $\wn$ is empty below this line, while Theorem
$\wt{{\rm  B}}$  says that the same is true for $\wnt$
except on $\mu=0$. On the  other hand, the Brill-Noether
number $\rn$ can be positive below the line, so  this is
not a sufficient condition for the non-emptiness of $\wn$.
This  phenomenon can be compared with the \lq\lq fractal
mountain range'' of Drezet  and Le Potier, which excludes
the existence of some stable bundles on ${\bf  P}^2$,
which should exist for purely dimensional reasons [DL].

\bigskip\bigskip \centerline{FIGURE 2}\bigskip

\renewcommand{\thesection}{\S\arabic{section}}
\section{Emptiness of Brill-Noether
loci}\renewcommand{\thesection}{\arabic{section}}

In this section, we assume  that $E$ has rank $n\geq 2$
and that either $E$ is stable and $\mu (E)\leq 1$ or   $E$
is   semistable and  $\mu (E)<1$. Our main purpose is to
prove the necessity  of the conditions in Theorems B and
$\wt{{\rm  B}}$ (see Theorem 3.3).

We begin with the following proposition, which will be
used many times in the  paper.

\bigskip

\noindent{\bf Proposition 3.1.} {\it Let $E$ be a stable
bundle of degree $d$,  $0\leq d \leq n$ (or a semistable
bundle with $0\leq d < n$), and $h^0(E) \geq k  > 0$. Let
$V$ be a subbundle of $E$ generated by $k$ independent
global sections   of $E$. Then $V$ is a trivial bundle of
rank $k$. } \bigskip

{\it Proof:} We have the exact sequence $$0\rightarrow
V\rightarrow E  \rightarrow F \rightarrow 0.$$  If $V$ is
non-trivial, there exists a section  $s\in h^0(E) $ such
that  $\deg D > 0, $ where $D$ is the divisor of zeros of
$s$. Then  $\deg L(D) > 0$ and $\mu (L(D))\geq 1.$  But
$L(D)$ is a subbundle of a  stable (resp. semistable)
bundle $E$ and $\mu (E) \leq1 (resp. <1)$. This is a
contradiction, so  $V\cong {\cal O}^k$.\epf

\bigskip

\noindent{\bf Remark 3.2.} i)  The above implies that for
any $E\in \wn $, \  $0\leq d \leq n$, $E$ can be presented
as an extension of the form   $$0\rightarrow {{\cal
O}^k}\rightarrow E \rightarrow F \rightarrow
0.\eqno(1)$$Similarly, every point of $\wnt$, $0\leq d<n$,
has a representative  $E$ which can be presented in this
form. \bigskip

ii) Note that, if $d\geq0$ and $E$ is stable, or $d>0$ and $E$ is semistable,
then $h^0(E^*)=0$. Except in the case $d=0$, $E$
semistable, we may therefore   assume that $h^0(F^*)=0$ in
the above sequence, as $F^*$ is a subbundle of  $E^*$.

\bigskip

\noindent{\bf Theorem 3.3. } {\it ${\cal W}^{k-1}_{n,d}$
is empty for $d>0$,  $n>d+(n-k)g$ and for $d=0$.
$\widetilde {\cal W}^{k-1}_{n,d}$ is empty for $d>0$,
$n>d+(n-k)g$  and for $d=0$, $k>n$.}

\bigskip
Theorem 3.3 has also been proved in the case $d>0$ by Anne
Maisani.\bigskip

{\it Proof.}:\  By Remark 3.2(i), every point of $\wnt$
can be represented by a  bundle $E$ of the form (1). It
follows at once that $\wnt$ is empty if $k>n$ or  if $k=n$
and $d>0$. Moreover, if $d=0$, (1) contradicts the
stability of $E$; so $\wn$ is empty. We can therefore
suppose that $k<n$ and $d>0$.

In this case, the extensions (1) are  classified by the
elements of the vector  space  $H=\bigoplus^k H^1(F^*)$,
i.e. by $k$-tuples $(e_1,\ldots,e_k)$ with  $e_i\in
H^1(F^*)$. Moreover two extensions are isomorphic if the
corresponding  points are in the same orbit of the natural
action of $GL(k)$ on $H$. Thus, if  $e_1,\ldots,e_k$ are
linearly dependent, we can suppose (using this action)
that  $e_k=0$; hence the extension has a partial splitting
to give $\co$ as a direct  summand of $E$, contradicting
the stability hypothesis.

Now, since $h^0(F^*)=0$ by Remark 3.2(ii), we have
$$h^1(F^*)=d+(n- k)(g-1).$$So  $e_1,\ldots,e_k$ are
necessarily linearly dependent if $k>d+(n- k)(g-1)$, or
equivalently $n>d+(n-k)g$.\epf \bigskip

 For future  convenience we finish this section with the
following proposition.

\bigskip

\noindent{\bf Proposition 3.4.} {\it Let $F$ be a fixed
bundle of rank $n-k$ and degree $d$ with $h^0(F^*)=0$.
Then, if $n\leq  d+(n-k)g$, the extensions $$0\rightarrow
{{\cal O}^k}\rightarrow E \rightarrow F  \rightarrow 0$$
with no trivial summands are classified up to automorphism
of  $\co^k$ by a variety of dimension $k(d+(n-k)g-n)$.}

\bigskip {\it Proof :} The extensions of this form are
classified by the  linearly independent $k$-tuples of
elements of $ H^1(F^*)$ modulo the linear  action of
$GL(k)$, in other words by the Grassmannian $Grass_k
(H^1(F^*))$.  Now\begin{eqnarray*}\dim Grass_k(H^1(F^*)) =
k(h^1(F^*)-k)&=& k(d+(n-k)(g-
1)-k)\\&=&k(d+(n-k)g-n).\end{eqnarray*}\epf

\bigskip

\renewcommand{\thesection}{\S\arabic{section}}
\section{Irreducibility}
\renewcommand{\thesection}{\arabic{section}}In this
section we shall use Proposition 3.1 and Remark  3.2 to
prove Theorems A and $\tilde{{\rm A}}$ when $0<d<n$. We
begin with a  lemma which is probably well known (and
certainly frequently assumed), but which  we could not
find in a suitable form in the literature (see [Ty,
Theorem 2.5.1] for a similar  result).

Let  ${\cal F} $ be a bounded set of non-stable bundles
of rank $n$ and degree  $d$. Then there exists a finite
number of families of bundles of rank $n$ over $X$,
parametrised by varieties $V_{\alpha}$, including
representatives of all bundles  in the given set (up to
isomorphism). For $v\in V_{\alpha}$, let $E_v$ denote the
corresponding bundle over $X$, and let $n_{\alpha,v}$
denote the dimension of the closure of the set $\{ w\in
V_{\alpha}|E_w\cong  E_v\}$ in $V_{\alpha}$.
Write$$m_{\alpha}=\min\{n_{\alpha,v}|v\in V_{\alpha}\},\
\ \ \ p=\max_{\alpha}\{\dim V_{\alpha}-m_{\alpha}\}.$$In
these circumstances, we  shall say that $\cf$ {\it depends
on at most $p$ parameters}.

\blem \label{311} {\it Any bounded set ${\cal F}$ of
non-stable bundles of rank $n$ depends  on at most
$n^2(g-1)$ parameters.} \elem \brem \label{32} i) Since
stable  bundles of rank $n$ and degree $d$ depend on
precisely $n^2(g-1)+1$ parameters,  this means that for
counting problems we can assume that the dimension of any
bounded family of vector bundles of rank $n$ is at most
$n^2(g-1)+1$.

ii) If $g=1$, Lemma 4.1  is not true. Actually there are \lq\lq more'' unstable
than stable bundles in this case (see [A]). \erem

{\it Proof of Lemma 4.1:} If $E$ is a non-stable vector
bundle of rank $n$  then  there exists a filtration $$
0=E_0\subset E_1\subset E_2 \subset ...\subset  E_r=E $$
with $E_i/E_{i-1} $ stable and $\mu (E_i/E_{i-1} ) \leq
\mu (E_{i- 1}/E_{i-2})$. For $E\in\cf$, the ranks and
degrees of the $E_i$ can take only  finitely many values,
so we can suppose these ranks and degrees are all fixed.

Let $$\rk(E_i/E_{i-1}) = n_i, \ \ \ \  \deg (Hom(E_j/E_
{j- 1},E_i/E_{i-1})) =  d_{j,i}$$ and let $\beta _i$ be
the minimum number of parameters on which the  set of
bundles which can occur as $E_i$ in the above filtration
depends. From  the exact sequence $$ 0\ra E_1\ra E_2\ra
E_2/E_1\ra 0 $$ we have that $$\beta _2  \leq n^2_1(g-1)+1
+n^2_2(g-1)+1$$if $h^1(Hom(E_2/E_1,E_1))=0$ for all
$E_1$,  $E_2/E_1$ and $$\beta _2 \leq n^2_1(g- 1)+1
+n^2_2(g-1)+1 +\max
\{h^1(Hom(E_2/E_1,E_1))\}-1$$otherwise, since $E_1$ and
$E_2/E_1$ are stable. In  the first case, clearly$$\beta
_2 \leq (n_1+n_2)^2(g-1).$$In the second, since
$Hom(E_2/E_1,E_1)$ is semistable and  $d_{2,1}\geq0$, we
have by Clifford's  theorem  $$h^0(Hom(E_2/E_1, E_1)) \leq
\frac{d_{2,1}}{2} + n_1n_2.\eqno(2)$$ So  by Riemann-Roch
$h^1(Hom(E_2/E_1,E_1)) \leq n_1n_2g - \frac{d_{2,1}}{2}$.

Therefore
$$ \begin{array}{lll} \beta _2& \leq &(n_1^2+n_2^2)(g-1)+1
+ n_1n_2g -  \frac{d_{2,1}}{2} \\ &=&(n_1+n_2)^2(g-1)+1
-n_1n_2(g-2) -\frac{d_{2,1}}{2}\\ &\leq
&(n_1+n_2)^2(g-1) \end{array} $$unless $g=2$ and
$d_{2,1}=0$. In the  exceptional case, the left-hand side
of (2) is $0$ unless $E_2/E_1\cong E_1$,  when it is $1$;
so the inequality can be improved unless $E_1$ and
$E_2/E_1$ are  isomorphic line bundles. But the extensions
of the required form in which $E_1$  and $E_2/E_1$ are
isomorphic line bundles depend
on$$2g-1\leq4(g-1)$$parameters.  This completes the proof
for $r=2$.

For $r\geq3$, we proceed by induction on $r$. The same
argument as above gives  $$\beta _r
\leq\beta_{r-1}+n_r^2(g-1)+1 +\max\{
h^1(Hom(E_r/E_{r-1},E_{r-1}))\}  -1$$(unless
$h^1(Hom(E_r/E_{r-1},E_{r-1}))$ is always zero, in which
case there  is a better estimate as above). Now
\begin{eqnarray*}h^1(Hom(E_r/E_{r-1},E_{r-
1}))&\leq&\sum_{i=1}^{r- 1}h^1(Hom(E_r/E_{r-1},E_i/E_{i-
1}))\\&\leq&\sum_{i=1}^{r-1}(n_in_rg-
\frac{d_{r,i}}{2})\end{eqnarray*} by  Clifford's Theorem
and Riemann-Roch.  So$$\beta _r \leq\beta_{r-1}+n_r^2(g-
1)+2\sum_{i=1}^{r-1}n_in_r(g-1),$$and the result  follows
from the inductive  hypothesis.\epf

\bigskip

We are now ready to prove Theorem $\tilde{{\rm A}}$ when
$0<d<n$.

\bth \label{33}  If $0<d<n$ and  $\wnt$ is non-empty, then
it is irreducible.  \eth

{\it Proof:} By Proposition 3.1 and Remark 3.2, any point
of $\wnt $ has a  representative $E$ of the form (1) with
$h^0(F^*)=0$. For fixed rank and degree,  the set $ \{ F |
h^0(F^*) =0\}$ is bounded. It follows by a standard
argument (due  originally to Serre, see for example [A,
Theorem 2]) that there is an irreducible  family which
includes representatives of all such bundles. The
condition  $h^0(F^*)=0$ defines an open subfamily,
parametrised by an irreducible variety  $Y$. The required
extensions are then parametrised by a projective bundle
over  $Y$, and those for which $E$ is semistable by an
open subset of the total space  of this bundle. This
subset is again irreducible and maps onto $\wnt$. So
$\wnt$  is irreducible.\epf

\bigskip Next we prove Theorem
A for $0<d<n$.

\bth \label{34}  If $0<d<n$ and $\wn$ is non-empty, then it is irreducible of
dimension $\rn$. Moreover $\Sing\wn =\wnn$. \eth

{\it Proof :}  Suppose $\wn $ is not empty. Since $\wn$ is an open subset of
$\wnt$, it is irreducible. From 1.1 we know that $\rn \leq $dim $\wn$.

Given $n$, $d$, $k$, let ${\cal S} $ be the set of all
possible extensions
$$0\ra \ok \ra  E \ra F \ra 0$$with  $h^0(F^*)=0$,
$\rk F=n-k$ and $\deg F=d$,  such that $E$  does not have
trivial summands. From Remark 4.2(i) and Propositions 3.1
and 3.4 we  obtain \begin{eqnarray*} \dim\wn & \leq &{\rm
the\ number\ of\ parameters\ on\ which\  }{\cal S}\ {\rm
depends}\\ &\leq&(n-k)^2(g- 1)+1 +k(d+(n-k)g-n)\\ &=&\rn
\end{eqnarray*}  Therefore, $\dim\wn = \rn$.

To see that $\Sing \wn = \wnn$, note first that, since
$\rn>\rho_{n,d}^k$,   $$\wn\neq\wnn.$$ Now let $E \in \wn
-\wnn$, so that $H^0(E) \cong H^0(\ok)$.  Since

$$ \begin{array}{lll} H^0(E)\otimes H^0(E^*\otimes K)&
\cong &H^0(\ok)\otimes  H^0(E^*\otimes K)\\ &\cong
&H^0(\ok \otimes E^*\otimes K)\\ &\hookrightarrow&
H^0(E\otimes E^*\otimes K), \end{array} $$ the Petri map
is injective. So, by  1.3, $\wn$ is smooth at $E$.
Since $\wnn\subset\Sing\wn$ by 1.3, we have
$\Sing\wn = \wnn$ as required.\epf

Theorems 4.3 and 4.4 complete the proof of Theorems A and  $\tilde{{\rm A}}$
when $0<d<n$. The cases $d=0$ and $d=n$ will be covered in
\S\S 7, 8.

We finish this section with \bcor If $(n-k,d)=1$ and
$\wn$ is non-empty, then  there is a dominant rational map
$g: Grass_k({\cal R}^1_p({\cal U}^*)) --\ra  \wn,$ where
${\cal U} $ is the universal bundle over $X\times {\cal
M}(n-k,d)$  and $p$ the projection to ${\cal M}(n- k,d)$.
\ecor

{\it Proof :} By Lemma 4.1 and the proof of Theorem 4.4,
the stable bundles $E$  constructed from non-stable $F$
belong to a proper subvariety of $\wn$. The  corollary now
follows from the proofs of Proposition 3.4 and Theorem
4.4. \epf

\bigskip
\renewcommand{\thesection}{\S\arabic{section}}\section{A
criterion for non-emptiness}\renewcommand{\thesection}
{\arabic{section}}

In this section we will give the setting that we need to
prove Theorems B and  $\tilde{\rm B}$ for $0<d<n$. More
precisely, we shall give a criterion for the
non-emptiness of $\wn$ by estimating the number of
conditions on an extension  (1) which are required for $E$
to be non-stable.

Assume that $0<d<n$ and let $F$ be a stable bundle of
rank $i$ and degree $d$.  Let $$ \xi : 0\ra \ok \ra E \ra
F \ra 0$$ be an extension of $F$ by $\ok$ such  that $E$
does not have trivial summands. By the proof of Theorem
3.3, such $\xi$  exist if and only if $k+i=n\leq d+ig.$

If $E$ is non-stable, then it has a stable quotient bundle
$H$ of rank $s<n$ and  degree $d'$ such that $$\mu (H)
\leq \mu (E).  \eqno(3)$$ This fits in the  following
diagram:

$$ \begin{array}{ccccccccc} 0&\ra &\ok&\ra &E&\ra &F&\ra
&0 \\  &&\da&&\da&&\da&&\\ 0&\ra &M&\ra &H&\ra &H_1&\ra &0
\\ &&\da&&\da&&\da&&\\  &&0&&0&&0&& \end{array}\eqno(4) $$
where $M$ is the image of $\ok \ra E\ra H$.  Note that
$M\not=0$; otherwise there would exist a
non-zero homomorphism $g:F\ra H$. Since  both are stable,
$\mu (F) \leq \mu (H)$. Since $\mu (E) < \mu (F)$, this
contradicts (3). Moreover, since $M$ has non-negative
degree and $H$ is stable, $\deg H\geq 0$

Since $M$ is generated by its global sections, it must be trivial. Otherwise
there would exist a section of $H$ generating a line bundle
of positive degree;  in conjunction with (3), this
contradicts the stability of $H$. For the same  reason,
$H_1$ must be torsion-free (and hence locally free).

One can now complete diagram (4) as follows

$$ \begin{array}{ccccccccc} &&0&&0&&0&&\\ &&\da&&\da&&\da
&&\\ 0&\ra &{\cal  O}^{n-s-l}&\ra &G&\ra &G_1^l&\ra &0 \\
&&\da&&\da&&\da&&\\ 0&\ra &\ok&\ra  &E^n&\ra &F^i&\ra &0
\\ &&\da&&\da&&\da&&\\ 0&\ra &{\cal O}^{s-m}&\ra &H^s&\ra
&H_1^m&\ra &0 \\ &&\da&&\da&&\da&&\\ &&0&&0&&0&&
\end{array}\eqno(5) $$(where  the superscripts denote the
ranks of the various bundles). Note that $m>0$;  otherwise
$H$ would be trivial  and $E$ would have a trivial
summand. Also  $l>0$; otherwise$$\deg H=\deg  H_1=\deg
F=d,$$contradicting (3).

Note that the existence of the top sequence in (5) implies
that the $k$-tuple of  elements of $H^1(F^*)$ defining
$\xi$ maps under the surjective homomorphism   $H^1(F^*)
\ra H^1(G_1)$ to a $k$-tuple of which at most  $(n-s-l) $
components   are linearly independent. Since, by
Riemann-Roch,$$h^1(G_1^*)\geq d-d'+l(g-  1),$$this rank
condition defines a subvariety $Z$ in  $\bigoplus^k
H^1(F^*)$   with $$ {\rm codim}\,Z \geq (s-m)(d-d'+lg
-n+s).\eqno(6)$$

On the other hand, the stability of $F$ implies that
every quotient bundle of   $H_1$ has slope greater than
every subbundle of $G_1$, and hence that
$h^0(H_1^*\otimes G_1)=0$. So$$h^1(H_1^*\otimes
G_1)=ld'-m(d-d')+lm(g-1).$$Since   $F$ varies in a bounded
set, so do $G_1$ and $H_1$ (note that $d-d'>0$ by (3));
so, by Lemma 4.1, the non-trivial extensions occuring in
the right-hand column   of (5) depend on at most

$ l^2(g-1) +1 +m^2(g-1) +1 + ld'-m(d-d')+lm(g-1) -1$

\hfill$=\dim{\cal M}(i,d)  + ld' -m(d-d')-lm(g-1)$\ \ \
\ (7)

\noindent  parameters.

\bprop \label{41} Let $0<d<n$. If $$(s-m)(d-d'+lg -n+s) >
ld'-m(d-d')-lm(g- 1)$$for all possible choices of $s$,
$d'$, $m$ and $l$, then $\wn $ is non-empty. \eprop

{\it Proof:} By (7), the general element $F$ of $\cm(i,d)$
admits families of  extensions $0\ra G_1\ra F\ra H_1\ra0$
as above depending on at most$$ld'-m(d-
d')-lm(g-1)$$parameters. If the inequality holds, it
follows from (6) that there  exists a non-empty open set
of extensions $\xi$  for which no diagram (5)  exists. If
this holds for all possible choices of $s$, $d'$, $m$ and
$l$ (of  which there are finitely many), then the general
extension $\xi$ must define a  stable bundle $E$.\epf

\bigskip

We will use Proposition 5.1 to prove Theorems B and
$\tilde{{\rm B}}$ for  $0<d<n$. In view of Theorem 3.3, it
is sufficient to show that the inequality of  Proposition
5.1 holds whenever the numerical conditions needed for (5)
to exist  hold. For convenience, we restate these
conditions now.

In the first place, (3) can be stated as $$sd-nd' \geq 0.
\eqno(a)$$ The  stability of $F$ implies that $$ (l+m)d'
-md >0. \eqno(b)$$         Since $H$ is  stable, we have
from the proof of Theorem 3.3 $$d' -s +mg \geq 0.
\eqno(c)$$ Finally the  inequality in Proposition 5.1 can
be written as $$m(n-s-l)-d'(l+s)+s(d+lg-  n+s)>0.
\eqno(d)$$

In the next section, we shall prove that $(a)$, $(b)$ and
$(c)$ imply $(d)$,  thus completing the proofs. \brem The
necessary condition $$n\leq d+(n-k)g $$ of  Theorems B and
$\tilde{\rm B}$ does not enter the calculation explicitly.
In  fact this inequality is a consequence of the
hypotheses of Proposition 6.1  below. \erem

\bigskip
\renewcommand{\thesection}{\S\arabic{section}}\section{Proof
of the inequality}\renewcommand{\thesection}{\arabic{section}} Our object in
this section is to prove
\bprop  Suppose $(a)$, $(b)$ and $(c)$ hold with $$0<d<n,\
\ 0<s\leq n-l\ \ {\rm and}\ \  l>0.$$ Then $(d)$ holds.
\eprop

It will be helpful for the proof to represent some of the data in a geometrical
form. We do this as follows:

\bigskip\bigskip \centerline{FIGURE 3}

\bigskip \bigskip

In this figure, we regard $n$, $d$, $s$ and $l$ as fixed
and $m$, $d'$ as  variables. The curve  $q$ is (a branch
of ) the hyperbola  $$(l+m)d'-md=0$$  defining the
inequality $(b)$ and has the form indicated since $l>0.$
The lines  $\ell _a $ and $\ell _c $ defining the
inequalities $(a)$ and  $(c)$ depend on  $s$, but the line
$\ell $ joining $C$ (the intersection of $\ell_a$ and
$\ell_c$) to $(0,0)$ has  equation $$(n-d)d' =dmg$$ which
is independent of $s$.  The shaded region is the region
where $(a),(b)$ and $(c)$ are all satisfied.
\blem
Suppose the hypotheses of Proposition 6.1 hold. Then
$$(n-d)s\geq n(n-d- lg).\eqno{(*)}$$\elem
{\it Proof:}
Note that, for any $s$, the line $\ell _c$   has slope $-g
<0.$ It follows that, if $(a),(b)$ and $(c)$ hold, then
the point  $C$ must lie above $q$.

Now $\ell $ meets $q$ at $(0,0)$ and the point $D$ with
coordinates  $$m=\frac{n-d-lg}{g} , \ \ \ \
d'=\frac{d(n-d-lg)}{n-d}.$$ Since $\ell $ has  positive
slope, the $m$-coordinate of $C$ must be at least as great
as that of  $D$, i.e. $$\frac{(n-d)s}{ng} \geq
\frac{n-d-lg}{g}.$$ Clearing denominators,  this gives
$(*)$. (Note that the coordinates of $D$ in this proof
could be  negative; this does  not affect the argument.)
\epf

{\it Proof of Proposition 6.1:}  The value of the LHS of
$(d)$  at $C$ is $$\frac{(n-d)s}{ng}(n-s-l)-
\frac{ds}{n}(l+s) +s(d+lg-n+s).$$ A simple calculation
shows that this is equal  to
$$\frac{s(g-1)}{ng}[s(n-d)-n(n-d-lg)+l(n-d)].$$ It follows
at once from Lemma  6.2 that this is positive.  In other
words, $C$ lies below the line defining the  inequality
$(d)$, which has non-negative slope. So the whole region
in which  $(a)$, $(b)$ and $(c)$ all hold also lies below
this line.

\epf

We are now ready to state\bth Theorems {\rm B} and
$\wt{{\rm  B}}$ hold for  $0<d<n$.\eth

{\it Proof:} This follows from Theorem 3.3 and
Propositions 5.1 and 6.1.\epf

\bigskip

\renewcommand{\thesection}{\S\arabic{section}}
\section{The case $\mu = 0 $}\renewcommand{\thesection}
{\arabic{section}}
\bth Theorems {\rm A} and {\rm B} hold for $d=0$.\eth

{\it Proof:} For
bundles of degree $0$, the existence of a  section
contradicts stability; so $\wn$ is always empty. This
gives Theorem B, and Theorem A holds trivially.\epf

\bigskip On the
other hand, we have \bth For $1\leq k\leq n,$ there
exists a bijective morphism $$\widetilde{\cm}(n-k,0)\ra
\wt{\cw }^{k- 1}_{n,0}.$$\eth {\it Proof:} If $E$ is a
semistable bundle of degree $0$ with  $k$ independent
sections, then by Proposition 3.1 we have an extension
$$0\ra  \ok \ra E\ra F\ra 0.$$ So $E$ is S-equivalent to
$\ok \oplus F$ for some  semistable bundle $F$ of rank
$n-k$ and degree $0$.

Hence the formula $[F] \ra [\ok \oplus F] $ defines a
bijection from  $\widetilde{\cm} (n-k,0)$ to $\wt{\cw
}^{k-1}_{n,0}.$ Since $\widetilde{\cm}  (n-k,0)$ is a
coarse moduli space, this is a morphism.\epf

\bth Theorems $\wt{{\rm  A}}$ and $\wt{{\rm  B}}$ hold for
$d=0$.\eth

{\it Proof:} By Theorem 3.3, $\wt{\cw}^{k-1}_{n,0} $ is
empty if $k>n$. The rest
of Theorem $\wt{{\rm  B}}$ now follows from Theorem 7.2, as
does Theorem  $\wt{{\rm  A}}$ when we recall that
$\widetilde{\cm} (n-k,0)$ is irreducible.\epf

 \brem Note that $$ \dim\widetilde{\cm} (n-k,0) =(n-k)^2
(g-1)+1 < \rho ^{k- 1}_{n,0}$$ if $n<(n-k)g$. This is
no contradiction since the points of
$\wt{\cw}^{k-1}_{n,0} $ correspond to S-equivalence
classes of bundles, not  isomorphism classes. \erem

\bigskip\renewcommand{\thesection}{\S\arabic{section}}\section{The
case $\mu = 1 $}\renewcommand{\thesection}
{\arabic{section}} In this final section we prove our
theorems for  the case $d=n$.

For stable bundles the key result is \bprop $\cw ^{n-2}_
{n,n} $ is non-empty.  \eprop {\it Proof:} Consider the
extensions $$0\ra \co ^{n-1} \ra E \ra F \ra  0,$$ where
$F$ is a line bundle of degree $n$. Since $n\leq n+g,$
there exist  extensions of this form for which $E$ has no
trivial summands. If $E$ is non-stable, we have as in \S
5 a diagram

$$ \begin{array}{ccccccccc} 0&\ra &\co ^{n-1}&\ra &E&\ra
&F&\ra &0 \\  &&\da&&\da&&\da&&\\ 0&\ra &M&\ra &H&\ra
&H_1&\ra &0 \\ &&\da&&\da&&\da&&\\  &&0&&0&&0&&
\end{array} $$ with $H$ stable and $\mu (H) \leq \mu (E)
.$ If $H_1$  is a line  bundle, then $H_1 \cong F$; so
$$\deg H=\deg F +\deg M \geq d$$ and  $\mu (H) > \mu (E)$,
which is a contradiction.

It follows that $H_1$ must be a torsion sheaf. If $\mu (H)
<1,$ this contradicts  the stability of $H$ just as in \S
5. However, if  $\mu (H) =1$, it is possible  for $H$ to
have a section with a zero. This can  happen only if
$H=\co (x) $ for  some $x \in X.$ Moreover, in this case,
we cannot have $M=\co (x),$ since $\co  (x) $ is not
generated by global sections, so our diagram must become
$$  \begin{array}{ccccccccc} &&0&&0&&0&&\\
&&\da&&\da&&\da&&\\ 0&\ra &{\cal O}^{n- 2}&\ra &G&\ra
&F(-x)&\ra &0 \\ &&\da&&\da&&\da&&\\ 0&\ra &\co ^{n-1}&\ra
&E&\ra  &F&\ra &0 \\ &&\da&&\da&&\da&&\\ 0&\ra &{\cal
O}&\ra &\co (x)&\ra &\co _x&\ra &0  \\ &&\da&&\da&&\da&&\\
&&0&&0&&0&& \end{array} $$

The existence of this diagram implies that the
$(n-1)$-tuple of elements of  $H^1(F^*)$ corresponding to
the extension $$0\ra \co ^{n-1} \ra E \ra F \ra 0$$  must
become dependent in $H^1(F(-x)^*)$. Now $$h^1(F(-x)^*) =
n+g-2 \geq n-1;$$  so this condition defines a subvariety
of $\bigoplus^{n-1} H^1(F^*)$ of  codimension  $g>1$. Since
$F(-x)$ depends on only one parameter, we can find an
extension for  which no such diagram exists.\epf

\bth Theorems {\rm A}  and {\rm B} hold for $d=n$.\eth

{\it Proof:} Proposition 3.1 remains true for stable
bundles when $d=n$. The  arguments of \S 4 therefore apply
to prove Theorem A in this case.

For Theorem B, it follows from Proposition 8.1 that $\cw
^{k-1}_{n,n}$ is non- empty for $k\leq n-1$. On the other
hand $\cw ^{n-1}_{n,n}$ is certainly empty,  since a
bundle with $n$ independent sections is either trivial or
has a section  with a zero; when $d=n$, either possibility
contradicts stability. Thus Theorem  B holds when
$d=n$.\epf

\bigskip We now turn to the semistable case. As in the case
$d=0$,  we obtain a result  which is interesting in its
own right. \bth Let $S^nX$ denote the $n$th symmetric
power of $X$. Then there exists a  bijective morphism
$S^nX \ra \wt{\cw}^{n- 1}_{n,n}.$ \eth {\it Proof:} Let
$E$ be a semistable bundle of rank and degree  $n$ with $n$
independent sections. Since $E\not\cong \co ^n,$ it must
possess a  section  with a zero. Semistability then gives
an extension $$0\ra \co (x) \ra E  \ra E' \ra 0,$$ where
$E'$ is semistable of rank  and degree $n-1$ and has  $n-
1$ independent sections. It follows by induction that $E$
is S-equivalent to a  bundle of the form $\co (x_1)\oplus
\dots \oplus \co (x_n).$ The existence of  the required
morphism follows from the universal properties of $S^nX$
and $\cm  (n,n)$.\epf

We need also \bprop For $k<n$, the point $[E] \in
\widetilde{\cw}^{k-1}_{n,n} $  determined by a semistable
bundle $E$ lies in the closure of $\cw ^{k-1}_{n,n}.$
\eprop {\it Proof:} For those bundles $E$ which can be
expressed as extensions  $$0\ra \ok \ra E\ra F\ra 0,$$ we
argue exactly as in Theorem 4.3.

The remaining bundles are those which possess a section
with a zero. We then  have an extension $$0\ra \co (x) \ra
E\ra F\ra 0,$$ so that $E$ is S-equivalent  to $\co (x)
\oplus F.$ We can suppose inductively that $[F]$ belongs
to the  closure of $\cw ^{k-2}_{n-1,n-1}$. (Note that, in
the case $n=2$, $F $ is a line  bundle, so this is
trivial. More generally, it is trivial whenever $k=1$,
since then $\cw ^{k-2}_{n-1,n-1}={\cm}(n-1,n-1)$.) It is
therefore sufficient to prove the proposition when
$E\cong \co (x) \oplus F$ and $F \in \cw ^{k-2}_{n-1,n-1}$.

For this, we consider extensions $$0\ra F \ra E'\ra \co
(x)\ra 0,$$ Note that we  have an inclusion $\co \subset
\co (x) $ and that this section of $\co (x)$  lifts to
$E'$ if and only if the pull-back of the extension by this
inclusion is  trivial. We therefore have a family of such
extensions parametrised by $$V=Ker  [H^1(\co (x)^* \otimes
F) \ra H^1(F)].$$ Now, since $F$ is stable with $\mu(F)
=1,$  $$h^1(\co (x)^* \otimes F) = (n-1)(g-1)$$and
$$h^1(F)=h^0(F) -(n-1) +(n-1)(g-1)  <(n-1)(g-1).$$ So
$\dim V\geq1$ and there exist non-trivial extensions
$$0\ra F  \ra E'\ra \co (x)\ra 0$$ with $h^0(E') \geq k$.

Now suppose $E'$ is such an extension, and that it
possesses a  section with a  zero. Since $F$ is stable,
this cannot be a section of $F$, so it maps to a  section
of $\co (x).$ The corresponding subbundle must map
isomorphically to  $\co (x),$ splitting the extension.

It follows that $V$ parametrises a family such that the
general member has  a  subbundle $\ok $ and therefore
defines a point in the closure of  $\cw ^{k- 1}_{n,n} $,
while the special member corresponding to $0 \in V $ is
$\co (x)  \oplus F.$ Hence $[\co (x) \oplus F]$ is in the
closure of $\cw ^{k-1}_{n,n} $  as required.\epf

We now have finally
\bth Theorems $\wt{{\rm  A}}$ and $\wt{{\rm  B}}$ hold
for $d=n$.\eth

{\it Proof:} $\widetilde{\cw} ^{k-1}_{n,n} $ is irreducible
for $k=n$  by Theorem 8.3 and for  $k<n$ by Theorem 8.2 and
Proposition 8.4. On the other hand, the proof of  Theorem
8.3 shows that no semistable bundle with $d=n$ can have
more than $n$  independent sections.\epf

\pagebreak\noindent {\small L. Brambila Paz\\*Departamento
de  Matematicas\\*UAM - Iztapalapa\\*Mexico, D. F.\\*C. P.
09340\\*Mexico\\*lebp@xanum.uam.mx \\[.15in] I. Grzegorczyk\\*Department of
Mathematics\\*University of Massachusetts\\*Dartmouth\\*MA
02747\\*U.S.A.\\*grze@cis.umassd.edu \\[.15in] P. E. Newstead\\*Department of
Pure Mathematics\\*The University of Liverpool\\*Liverpool\\*L69
3BX\\*England\\*newstead@liverpool.ac.uk}

\end{document}